\documentclass[preprint,12pt]{elsarticle}

\usepackage{amssymb}
\usepackage{amsmath}

\usepackage{subcaption}
\usepackage{algorithm}
\usepackage{algorithmic}

\begin{document}

\begin{frontmatter}



\title{Phase Retrieval Based on DC and DnCNN}


\author[label1]{Xueming Li} 
\author[label1]{Bing Guo}

\affiliation[label1]{organization={School of Mathematics and Statistics},
            addressline={Jishou University}, 
            city={Jishou},
            postcode={416000}, 
            state={Hunan Province},
            country={China}}

\begin{abstract}
This paper investigates noise-robust phase retrieval by enhancing the prDeep architecture with difference of convex functions (DC) and DnCNN-based denoising regularization. This research introduces two novel algorithms, prDeep-DC and prDeep-L2, which demonstrably achieve excellent quantitative and visual performance, as confirmed by extensive numerical experiments.
\end{abstract}

%
%
%



\end{frontmatter}



\section{Introduction}

In numerous practical scenarios, signals invariably encounter degradation due to noise, distortion, or various forms of interference. This corruption significantly compromises the integrity of the phase information. Phase retrieval endeavors to reconstruct the original signal,   $\mathbb{R}^n$ or $\mathbb{C}^n$ , from a set of measurements, $\mathbf{b}=\left(b_1,\cdots,b_n\right)$, which are typically of the form
\begin{equation}
	\mathbf{b}=\vert A\mathbf{x}\vert +w.
\end{equation}
where A is the measurement matrix and $w$ represents noise.

The phase retrieval problem is inherently defined as an inverse problem: the recovery of a signal from the amplitude of the output of a nonlinear system. Such inverse problems are generally formulated as optimization problems of the following form
\begin{equation}
	\min_\mathbf{x}  \ell \left(f\left(\mathbf{x}\right),\mathbf{b}\right)+\lambda R\left(\mathbf{x}\right),
\end{equation}
where $f\left( \cdot \right)$ is the forward operator, $\ell \left(f\left(\mathbf{x}\right),\mathbf{b}\right)$ denotes the data fidelity term,
$R\left(\mathbf{x}\right)$ denoting the prior information of $\mathbf{x}$ can be used to solve the optimization problem as the regularization term of the objective function, and $\lambda$ denotes the degree of penalty to $R\left(\mathbf{x}\right)$. For the fourier phase recovery problem, it can be written in the following form
\begin{equation}
	\min_{\mathbf{x}} \, \Vert \vert A\mathbf{x}\vert-\mathbf{b} \Vert ^2_2+\lambda R\left(\mathbf{x}\right),
\end{equation}
here $A$ is the discrete fourier transform matrix.

Phase retrieval algorithms originated in the 1970s with the Gerchberg-Saxton (GS) algorithm \cite{gerchberg1972practical}. Fienup’s Hybrid Input-Output (HIO) algorithm \cite{fienup1984reconstruction} later improved upon GS and remains a widely used convex optimization-based approach. A systematic understanding of the connection between phase retrieval and convex optimization was established in \cite{2002Phase}. Simultaneously, phase retrieval has shown strong ties to compressive sensing, with research in \cite{moravec2007compressive} and \cite{schniter2014compressive} focusing on sparse phase retrieval by integrating compressive sensing concepts.

In recent years, significant advancements in tackling the phase retrieval problem have been propelled by the adaptation of both the Alternating Direction Method of Multipliers (ADMM) and neural network architectures. Specifically, the ADMM algorithm has been employed to address a variety of phase retrieval challenges, as demonstrated in \cite{wen2012alternating}. While ADMM exhibits a notable degree of robustness with respect to the selection of relaxation parameters, a considerable challenge persists in determining appropriate multiple parameters for optimal performance.


The rapid progress in deep learning \cite{zuo2022deep} has led to the widespread adoption of deep neural networks in image processing. In phase retrieval, \cite{tayal2020unlocking} utilizes them to train a mapping from observed signals to their expected reconstructions. By breaking phase retrieval symmetries with deep neural networks, phase information can be recovered from amplitude data. A key limitation, however, is the method’s suitability only for images related to the training sets. To overcome such data-dependent limitations, the Regularization by Denoising (RED) was introduced in \cite{romano2017little}. RED is defined as follows
\begin{equation}
	R\left(\mathbf{x}\right)=\frac{\lambda}{2}\mathbf{x}^T\left(\mathbf{x}-D\left(\mathbf{x}\right)\right) \label{RED},
\end{equation}
where $D\left(x\right)$ is an arbitrary denoiser.


RED can integrate any image denoising algorithm for image recovery, endowing it with excellent flexibility and applicability. In \cite{metzler2018prdeep}, prDeep, the robust phase retrieval method with a flexible deep network, was proposed by applying RED to noisy phase retrieval tasks. Simulation results demonstrate its noise robustness, yet it still exhibits insufficient resilience to complex noise. Since the loss function $\ell \left(f\left(\mathbf{x}\right),\mathbf{b}\right)$ is
is typically nonconvex and can be expressed as the difference of convex functions, the Difference of Convex Algorithm (DCA) is adopted for phase retrieval. In \cite{2015PhaseLiftOff}, PhaseLiftOff was developed as an improved version of PhaseLift\cite{2011PhaseLift}, a widely used phase retrieval algorithm.

In our work, we have made two contributions. The first one is that we apply DCA to replace the first step of the proximal gradient method and propose an improved alogrithm, which is denoted by prDeep-DC algorithm. The second contribution is that we study the 2-norm version of the regularization by denoising and analyze the relationship between the 2-norm regularization and regularization by denoising in RED\cite{romano2017little}, and propose prDeep-L2 algorithm. Finally, we verified the superiority of our improved method both in quantitative values and visual effects through numerical experiments. 


The remainder of this paper is structured as follows. Section 2 provides background on relevant concepts, including the image denoising engine, DC (difference of convex functions) formulation, the forward-backward splitting (FBS) algorithm, and its enhanced variants. In Section 3, we present our main contributions, which involve the application of DCA and L2 regularization via denoising, along with the corresponding improved algorithms. Section 4 details the experimental setup and analyzes the numerical results. Finally, Section 5 concludes the paper with a summary of our work.

\section{Preliminary knowledge}
\subsection{Image denoising engine}
Image denoising is an important task in digital image processing, which aims to eliminate or reduce noise from the image and improve the quality and clarity of the image. The goal of image denoising is to recover a clean image $\mathbf{x}$ from a noisy observation image $\mathbf{y}$. Its model can be written as following
\begin{equation}
	\mathbf{y} = \mathbf{x} + \mathbf{v} ,
\end{equation}
where $\mathbf{x}$ is the unknown clean image, $\mathbf{y}$ is the observed image with noise, $\mathbf{v}$ is generally additive white Gaussian noise or Poisson noise.  It can be seen from this model that image denoising is an inverse problem, then it remains a challenge and open problem. 

The classical denoising methods depend on the various filters, which can be classfied into spatial domain filter and tranform domain filter. The spatial domain filter can be divided into linear and nonliear filters, and tranform domain filter includes fourier transform, Gabor transform, wavelet transform, curvelet transform and so on. A large number of efficient denoising algorithms have been proposed over the past decades of years . In fact, current denoising algorithms mainly combine more than one filter with related knowledge. Here, we will introduce BM3D and DnCNN, which will be used in our later numerical experiments.

\begin{figure}[t]
	\centering
	\renewcommand{\figurename}{Figure}
	\includegraphics[width=0.9\textwidth]{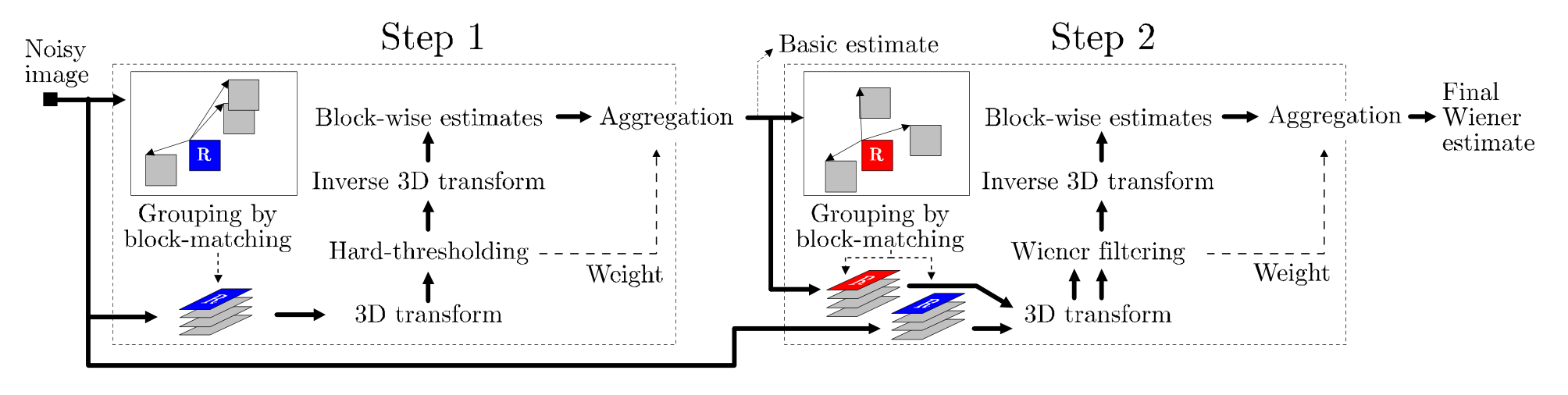}
	\hfill
	\caption{BM3D.}
	\label{BM3D}
\end{figure}

 BM3D(Block-Matching and 3D Filtering)\cite{dabov2007image} is a denoising algorithm based on block matching and 3D filtering. It divides the image into overlapping blocks and constructs a block matching table by looking for similar blocks, which are shown in Figure \ref{BM3D}. It requires several steps such as block matching, clustering and 3D filtering, and each step involves a lot of computation operations. Therefore, the calculation time and memory consumption of BM3D algorithm are large. and its performance depends on  some key parameters in BM3D algorithm such as the block size, search window size, threshold and so on.

\begin{figure}[t]
	\centering
	\renewcommand{\figurename}{Figure}
	\includegraphics[width=0.9\textwidth]{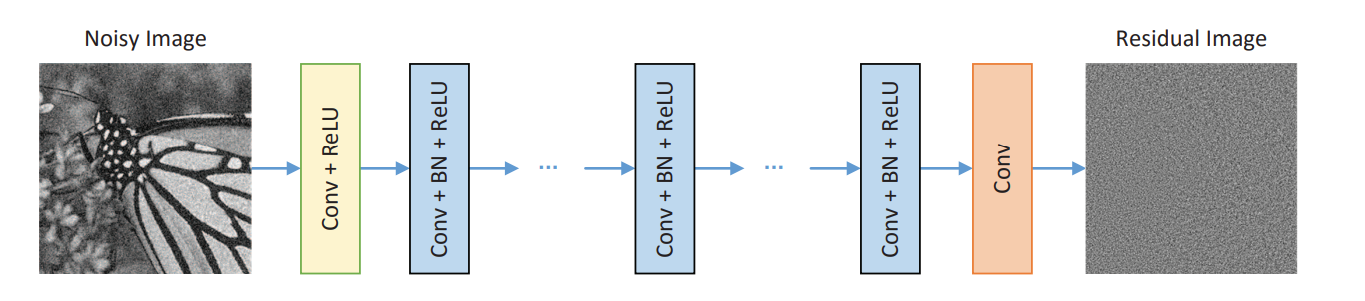}
	\hfill
	\caption{DnCNN.}
	\label{DnCNN}
\end{figure}

DnCNN (Denoising Convolutional Neural Networks) \cite{zhang2017beyond} is a learning-based denoising algorithm. By training on a large dataset of noisy-clean image pairs, it learns a mapping from noisy inputs to clean images, thereby achieving effective denoising. Although DnCNN requires substantial labeled data during training, it demonstrates strong performance in practical applications, capable of removing various types of noise while accurately restoring image details and textures. Additionally, DnCNN incorporates residual connections to facilitate faster model convergence. The input to DnCNN is a noisy image, and the output is a denoised image of the same size and channel number. The network comprises multiple convolutional layers, each consisting of convolution, an activation function, and batch normalization. The network structure diagram of DnCNN is showed in Figure \ref{DnCNN}.

The Regularization by Denoising (RED) framework is a widely adopted image reconstruction method. It introduces a regularization term to constrain the smoothness or sparsity of the solution, thereby removing noise. The fundamental principle is to balance data fidelity and regularization terms for optimal denoising. Like other plug-and-play algorithms, RED can utilize any denoiser for any image inverse problem. However, RED distinguishes itself by minimizing an explicit loss function with a denoiser, ensuring global optimality, unlike methods that implicitly minimize a loss. This transforms regularization into a denoising problem, leveraging powerful denoisers to restore clean signals, making RED highly robust to noise. Unlike traditional methods, RED requires no prior information or strict noise modeling; it adaptively adjusts parameters and recovers optimal solutions by learning noise statistics from the data.

The regularization term defined by regularization by denoising framework is 
\begin{equation}
	R\left( \mathbf{x} \right)=\frac{\lambda}{2} \mathbf{x} ^T\left( \mathbf{x} -D\left( \mathbf{x} \right)\right) , \label{RED}
\end{equation}
where $ D\left( \mathbf{x} \right) $ is an arbitrary denoiser and $\lambda$ is a regularization parameter. We use the traditional denoiser BM3D and the convolutional neural network-based DnCNN, and then apply the denoise regularization framework to solve the noisy phase retrieval in our experiments. 

\begin{figure}[t]
	\centering
	\begin{subfigure}{0.13\linewidth}
		\centering
		\includegraphics[width=\linewidth]{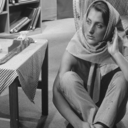}
		\caption*{Barbara}
	\end{subfigure}
	\centering
	\begin{subfigure}{0.13\linewidth}
		\centering
		\includegraphics[width=\linewidth]{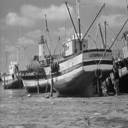}
		\caption*{Boat}
		
	\end{subfigure}
	\centering
	\begin{subfigure}{0.13\linewidth}
		\centering
		\includegraphics[width=\linewidth]{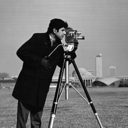}
		\caption*{Cameraman}
		
	\end{subfigure}
	\centering
	\begin{subfigure}{0.13\linewidth}
		\centering
		\includegraphics[width=\linewidth]{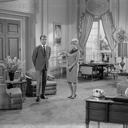}
		\caption*{Couple}
		
	\end{subfigure}
	\centering
	\begin{subfigure}{0.13\linewidth}
		\centering
		\includegraphics[width=\linewidth]{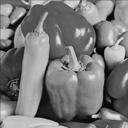}
		\caption*{Peppers}
		
	\end{subfigure}
	\centering
	\begin{subfigure}{0.13\linewidth}
		\centering
		\includegraphics[width=\linewidth]{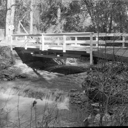}
		\caption*{Bridge}
		
	\end{subfigure}
	\caption{Test dataset NT-6 which are widely used in phase retrieval. }
	\label{NT-6}
\end{figure}

\begin{figure}[t]
	\centering
	\begin{subfigure}{0.13\linewidth}
		\centering
		\includegraphics[width=\linewidth]{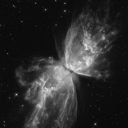}
		\caption*{Butterfly}
	\end{subfigure}
	\centering
	\begin{subfigure}{0.13\linewidth}
		\centering
		\includegraphics[width=\linewidth]{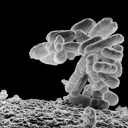}
		\caption*{Ecoli}
		
	\end{subfigure}
	\centering
	\begin{subfigure}{0.13\linewidth}
		\centering
		\includegraphics[width=\linewidth]{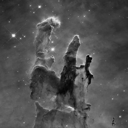}
		\caption*{Pillarsof} 
		
	\end{subfigure}
	\centering
	\begin{subfigure}{0.13\linewidth}
		\centering
		\includegraphics[width=\linewidth]{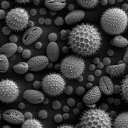}
		\caption*{Pollen}
		
	\end{subfigure}
	\centering
	\begin{subfigure}{0.13\linewidth}
		\centering
		\includegraphics[width=\linewidth]{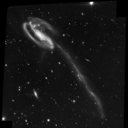}
		\caption*{Tadpole} 
		
	\end{subfigure}
	\centering
	\begin{subfigure}{0.13\linewidth}
		\centering
		\includegraphics[width=\linewidth]{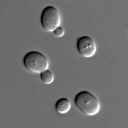}
		\caption*{Yeast}
		
	\end{subfigure}
	\caption{Test dataset UNT-6 which are widely used in phase retrieval.}
	\label{UNT-6}
\end{figure}

\subsection{DCA}
The following programming is called DC programming such that
\begin{equation}
		\min_{ \mathbf{x} } \,F\left( \mathbf{x} \right):= F_1\left( \mathbf{x} \right)-F_2\left( \mathbf{x} \right) ,
\end{equation}
where  $F_1\left( \mathbf{x} \right)$ and $F_2\left( \mathbf{x} \right)$
 are convex functions. This makes $F\left(\mathbf{x}\right)$
a DC (Difference of Convex) function, and  $F_1\left( \mathbf{x} \right)- F_2\left( \mathbf{x} \right)$ is its DC decomposition. As most nonconvex functions can be represented by DC functions, DC programming plays a crucial role in nonconvex optimization. The DCA was first introduced for DC programming in \cite{tao1986algorithms}. The core idea of DCA is to decompose the original nonconvex optimization problem into a series of subproblems, each being the difference of convex functions. DCA gradually approaches the optimal solution of the original problem by repeatedly solving these subproblems. The basic DCA scheme is presented in Algorithm 1.

The first step of DCA is an optimization problem in $\mathbf{x}$, solvable by algorithms like general gradient descent. DC programming and DCA have undergone substantial improvements and are now used to solve a variety of difficult nonconvex programs across fields like machine learning, signal processing, wireless communication, and image processing. For example, DCA is applied to train Support Vector Machines and logistic regression models in machine learning. In fact, the DC function has infinitely many DC decompositions which have great influence on the stability, rate of convergence and global solution of DCA. Therefore, it is important to choose a suitable DC decomposition.

\begin{algorithm}[H]
	\caption{DCA(DC Algorithm)}
	\label{DCA}
	\begin{algorithmic}[1]
		\REQUIRE Let $\mathbf{x}_0$ be any initial point and set k :=0.
		\STATE Select $\mathbf{u}_k \in \partial F_2(\mathbf{x}_k)$ and solve the strongly convex optimization problem$$ \min_{\mathbf{x}} \, F_1\left( \mathbf{x} \right)-\textless \mathbf{u}_k,\mathbf{x} \textgreater  $$ 
		to obtain its unique solution $\mathbf{y}_k$.
		\IF{ $\mathbf{y}_k=\mathbf{x}_k$} 
		\STATE STOP and RETURN $\mathbf{x}_k$.
		\ELSE 
		\STATE Set $\mathbf{x}_{k+1} :=\mathbf{y}_k $, set $k :=k+1$, and go to Step 1.
		\ENDIF 
	\end{algorithmic}
\end{algorithm}

\subsection{FBS algorithm and its improvement}
We consider the following optimization problem
\begin{equation}
	\min_{\mathbf{x}} \, f\left( \mathbf{x} \right)+g\left( \mathbf{x} \right) , \label{f1}
\end{equation}
where $ \mathbf{x} \in R^N$, $f \left( \cdot \right) $ is convex and differentiable function, and g is any convex function. Since the function $g \left( \cdot \right) $ is often nonsmooth, the standard gradient descent algorithm cannot be directly applied. Instead, we employ the Forward-Backward Splitting (FBS) algorithm, also known as proximal gradient descent. The proximal operator for the function $g \left( \cdot \right) $ is calculated as

\begin{equation}
	\rm{prox}_g\left( \mathbf{z},\tau \right) = {\arg\min_{\mathbf{x}}} \, \tau g\left( \mathbf{x} \right)+\frac{1}{2}\left\|\mathbf{x}-\mathbf{z}\right\|_2^2 , \label{prox_g}
\end{equation}
where $\tau$ is the step size.

%
By decomposing the problem (\ref{f1}) into forward and backward sub-problems, the FBS algorithm enhances solution feasibility and simplifies problem handling. This iterative approach gradually refines variable values in each iteration, approaching an approximated solution with improved results. The key advantage of FBS is its ability to effectively address combinations of differentiable and nonsmooth/non-differentiable parts in the objective function, as well as problems with separable additive structures, common in many applications. For convex problems, it converges to the global minimum. The specific FBS algorithm is outlined in Algorithm 
, making it suitable for a variety of optimization problems in signal processing, machine learning, and image processing.
\begin{algorithm}[H]
	\caption{FBS(Forward-Backward Splitting)\cite{goldstein2014field}}
	\label{FBS algorithm}
	\begin{algorithmic}[1]
		\REQUIRE Let $ \mathbf{x}_0 $ be an initial point, set k :=0 and set $\tau \textless 1/L\left( \nabla f \right)$.
		\FOR{$ k=1,2,3,\ldots ,$MAXITER}
		\STATE $\hat{\mathbf{x}}_k=\mathbf{x}_{k-1}-\tau \nabla f(\mathbf{x}_{k-1})$.
		\STATE $ \mathbf{x}_k=\rm{prox}_g\left(\hat{\mathbf{x}}_k,\tau\right)$.
		\ENDFOR
		\STATE Output current point $\mathbf{x}_{MAXITER}$ as the approximate optimal solution.
	\end{algorithmic}
\end{algorithm}

\begin{algorithm}[H]
	\caption{FISTA \cite{beck2009fast}} 
	\label{FISTA algorithm}
	\begin{algorithmic}[1]
		\REQUIRE $\mathbf{y}_1 = \mathbf{x}_0 \in R^n, \alpha_1=1, \tau \textless 1/L\left(\nabla f \right)$.
		\FOR{$ k=1,2,3,\ldots ,$MAXITER}
		\STATE $ \mathbf{x}_k = \rm{prox}_g\left(\mathbf{y}_k - \tau \nabla f\left( \mathbf{y}_k \right)\right)$.
		\STATE $ \alpha_{k+1} = \left( 1+\sqrt{1+4{\left(\alpha_{k}\right)}^2} \right)/2 $.
		\STATE $ \mathbf{y}_{k+1} = \mathbf{x}_k + \frac{\alpha_k -1}{\alpha_{k+1}}\left(\mathbf{x}_k - \mathbf{x}_{k-1}\right) $.
		\ENDFOR
		\STATE Output current point $\mathbf{x}_{MAXITER}$ as the approximate optimal solution.
	\end{algorithmic}
\end{algorithm}


Building upon the convergence of the FBS algorithm, numerous scholars have introduced improvements. The Fast Iterative Shrinkage-Thresholding Algorithm (FISTA) \cite{beck2009fast} is a key example, designed to accelerate FBS. FISTA enhances convergence by introducing an additional acceleration parameter, $\alpha$,
and employs a prediction-correction method. This method utilizes previous and current estimates to achieve faster convergence. The specific FISTA algorithm is detailed in Algorithm \ref{FISTA algorithm}.

In Algorithm \ref{FISTA algorithm}, $ L\left(\nabla f \right) $ is the Lipschitz constant of $\nabla f $, that is, it needs to satisfy
\begin{equation} 
	\left\| \nabla f\left( \mathbf{x}_1 \right) - \nabla f\left( \mathbf{x}_2 \right) \right\| \leq L\left\| \mathbf{x}_1 - \mathbf{x}_2 \right\|  \qquad   \forall \mathbf{x}_1, \, \mathbf{x}_2 .
\end{equation}
The algorithm converges when the step size is $\tau \textless 1/L\left(\nabla f \right)$. The second step of algorithm \ref{FISTA algorithm} is the original FBS algorithm, and the third step is to accelerate the update process of parameter $ \alpha$, where $ \alpha$ keeps increasing with the iteration process, and when the objective function increases, the restart method is taken, that is, $ \alpha = 1$. The fourth step is a prediction
process, which can accelerate the convergence speed of the algorithm.

When the properties of $ \nabla f$ are unknown, choosing a suitable step size is difficult. Therefore, the backtracking line search method is used to ensure convergence. This method checks the line search condition after each FBS iteration and reduces the step size if the condition is not met, repeating the FBS iteration until the condition is satisfied.

Although the line search method makes the objective function monotonously decline, but it is relatively conservative. Then non-monotone line search is proposed and it allows
the objective function to increase within a certain range. The non-monotone line search condition is looser than the monotone line search condition, which means that there are fewer backtracking steps, and then the running time will be reduced accordingly.

Let $M$ be a line search parameter which is a positive integer, and then define
\begin{equation}
	\hat{f}_{k}=\max\left \{ f_{k-1}, f_{k-2},\ldots, f_{k-min\left \{M,k  \right \} }   \right \} ,
\end{equation}
\begin{equation}
	f\left( \mathbf{x}_{k+1}\right) < \hat{f}_{k} + \mathfrak{Re} \big\langle \mathbf{x}_{k+1} - \mathbf{x}_{k},\nabla f\left(\mathbf{x}_{k}\right) \big\rangle + \frac{1}{2\tau_{k}}\left\|\mathbf{x}_{k+1} - \mathbf{x}_{k}\right\|^2 . \label{Line Search condition} 
\end{equation}

If the backtracking condition is violated (\ref{Line Search condition}), then the step size is reduced until it is
established, and the process of non-monotone line search is in algorithm \ref{Non-Monotone}. Moreover, the
FBS algorithm using non-monotone line search is convergent\cite{goldstein2014field}.

\begin{algorithm}[H]
	\caption{Non-Monotone Line Search}
	\label{Non-Monotone}
	\begin{algorithmic}[1]
		\WHILE{$ \mathbf{x}_{k} $ and $ \mathbf{x}_{k+1}$ violate condition(\ref{Line Search condition})}
		\STATE $\tau_{k} \gets \tau_{k}/2$.
		\STATE $ \mathbf{x}_{k+1} \gets \rm{prox}_g\left(\mathbf{x}_{k}-\tau_{k}\nabla f\left(\mathbf{x}_{k} \right),\tau_{k} \right) $.
		\ENDWHILE
	\end{algorithmic}
\end{algorithm}

In \cite{goldstein2014field}, the fast adaptive shrinkage-thresholding algorithm denoted by FASTA is proposed. Comparing with FISTA, it has adaptive step size to accelerate the solving process, which can make our problem converge quickly.

\section{Our work}
In this paper, we consider the following form of phase retrieval problem
\begin{equation}
	\min_{\mathbf{x}} \, \Vert \vert \mathrm{\mathbf{A}} \mathbf{x}\vert-\mathbf{b}\Vert_2^2 \\
	=\min_{\mathbf{x}} \, \sum_{i=1}^{n} \left(\vert \langle \mathrm{\mathbf{A}}_i,\mathbf{x} \rangle \vert -\mathbf{b}_i\right)^2 ,
\end{equation}
where $\mathbf{x} \in R^n$, $ \mathrm{\mathbf{A}}_i $ is ith row of $ \mathrm{\mathbf{A}} $. Let
\begin{equation}
	F(\mathbf{x})=\sum_{i=1}^{n} \left(\vert \langle \mathrm{\mathbf{A}}_i,\mathbf{x} \rangle \vert -\mathbf{b}_i\right)^2 ,
\end{equation}
be the objective function. Because $ F \left( \mathbf{x} \right) $ is a nonconvex function, we can write it in the following form
\begin{equation}
	F (\mathbf{x})=\sum_{i=1}^{n} \left(\vert \langle \mathrm{\mathbf{A}}_i,\mathbf{x} \rangle \vert -\mathbf{b}_i\right)^2 =F_1\left(\mathbf{x}\right)-F_2 \left(\mathbf{x}\right) ,
\end{equation}
where $ F_1\left(\mathbf{x}\right)=\sum_{i=1}^{n} \left(\vert \langle \mathrm{\mathbf{A}}_i,\mathbf{x} \rangle \vert ^2+ \mathbf{b}_i^2\right)$, $ F_2\left(\mathbf{x}\right)=\sum_{i=1}^{n} 2\mathbf{b}_i\vert \langle \mathrm{\mathbf{A}}_i,\mathbf{x} \rangle \vert $. Because $ F_1\left( \mathbf{x} \right)$ and $ F_2\left( \mathbf{x} \right)$ are convex functions, so we can consider the original objective function as the difference of two convex functions and use DCA to solve this optimization problem. For any given $ \mathbf{x}_{k}$, we can get the following iterative process that
\begin{equation}
	\mathbf{x}_{k+1} = {\arg\min_{\mathbf{x}}} \, F_1\left( \mathbf{x} \right)-\nabla F_2(\mathbf{x}_{k})^T(\mathbf{x}-\mathbf{x}_{k}) .
\end{equation}

We study the phase retrieval of the image contaminated by noise, and use RED(\ref{RED}). Then the model is the following 
\begin{equation}
	{\arg\min_{\mathbf{x}}} \, \frac{1}{2{\sigma}^2} \left\| \mid \mathrm{\mathbf{A}} \mathbf{x} \mid - \mathbf{b}\right\|_2^2+\frac{\lambda}{2} \mathbf{x}^T\left( \mathbf{x} - D\left( \mathbf{x} \right)\right) , \label{the original phase retrieval}
\end{equation}
where $ \frac{\lambda}{2} \mathbf{x}^T\left(\mathbf{x} -D\left( \mathbf{x} \right)\right) $ is regularization by denoising, $\lambda$ is a regularization parameter, $\sigma$ is the
standard deviation of noise, and $ D\left( \mathbf{x} \right) $ is the denoiser. 

Based on DC decomposition, we get the following optimization problem
\begin{equation}
	{\arg\min_{\mathbf{x}}} \, \frac{1}{2{\sigma}^2} F_1\left( \mathbf{x} \right)-\frac{1}{2{\sigma}^2} F_2\left( \mathbf{x} \right)+\frac{\lambda}{2} \mathbf{x}^T\left(\mathbf{x} - D \left(\mathbf{x} \right)\right) \label{the phase retrieval with DC} .
\end{equation}
and use the FASTA solver, which is an efficient and convenient solver for this optimization problem.  Let $f\left( \mathbf{x} \right)=\frac{1}{2{\sigma}^2} F_1\left(  \mathbf{x} \right)-\frac{1}{2{\sigma}^2} F_2\left( \mathbf{x} \right)$ and $g( \mathbf{x} )=\frac{\lambda}{2} \mathbf{x}^T \left(\mathbf{x}-D\left(\mathbf{x}\right) \right)$. Because $f\left(\mathbf{x} \right)$ here is not differentiable, we take a subgradient $ \frac{1}{\sigma^2}\left(\mathrm{\mathbf{A}} \mathbf{x}-\mathbf{b} \circ \frac{\mathrm{\mathbf{A}} \mathbf{x}}{\mid \mathrm{\mathbf{A}} \mathbf{x} \mid}\right) $, where the symbol $ \circ $ denotes multiplication of the corresponding element. For the second step of Algorithm \ref{prDeep-DC}, the DCA and gradient descent method are introduced to obtain $\hat{\mathbf{x}}_k$, refering to Algorithm
\ref{prDeep-DC} for details. The proximal operator can be used for the third step. It is the following iterative process that
\begin{equation}
	\mathbf{x}_0 = \hat{\mathbf{x}}_k, \; \mathbf{x}_j = \frac{1}{1 + \lambda \tau} \left( \hat{\mathbf{x}}_k + \lambda \tau D\left( \mathbf{x}_{j-1} \right) \right) \; \; j=1, \cdots,T ,
\end{equation}
where $T$ is iteration number. By the above iteration processing, we have
\begin{equation}
	\mathbf{x}_k = \mathbf{x}_T .
\end{equation}

It is found that only one iteration can show good results, which means that
\begin{equation}
	\mathbf{x}_k = \frac{1}{1+\lambda \tau} \left( \hat{\mathbf{x}}_k + \lambda \tau D\left( \hat{\mathbf{x}}_k \right) \right) .
\end{equation}

\begin{algorithm}[H]
	\caption{prDeep}
	\label{prDeep}
	\begin{algorithmic}[1]
		\REQUIRE $\mathbf{y}_{1} = \mathbf{x}_{0} \in R^N, \alpha_{1}=1, \tau \textless 1/L\left(\nabla f \right)$
		\FOR{$ k=1,2,3,\ldots ,$MAXITER}
		\STATE $ \hat{\mathbf{x}}_{k}= \mathbf{y}_{k}-\frac{\tau}{2{\sigma}^2}\nabla F\left(\mathbf{y}_{k}\right) $
		\STATE $ 
		\mathbf{x}_{k} = {\rm{prox}_g}(\hat{\mathbf{x}}_{k},\tau)={\arg\mathop{\min} \limits_{\mathbf{x}}} \, \frac{\lambda \tau}{2} \mathbf{x}^T\left(\mathbf{x}-D\left(\mathbf{x}\right)\right)+\frac{1}{2}\left\| \mathbf{x}-\hat{\mathbf{x}}_{k}\right\|_2^2 $.
		\STATE $ \alpha_{k+1} = \left( 1+\sqrt{1+4{\left(\alpha_{k}\right)}^2} \right)/2 $
		\STATE $ \mathbf{y}_{k+1} = \mathbf{x}_{k} + \frac{\alpha_{k}-1}{\alpha_{k+1}}\left(\mathbf{x}_{k} - \mathbf{x}_{k-1}\right) $
		\ENDFOR
		\STATE Output current point $\mathbf{x}_{MAXITER}$ as the approximate optimal solution.
	\end{algorithmic}
\end{algorithm}

\begin{algorithm}[!h]
	\caption{prDeep-DC}
	\label{prDeep-DC}
	\begin{algorithmic}[1]
		\REQUIRE $\mathbf{y}_{1} = \mathbf{x}_{0} \in R^N, \alpha_{1}=1, \tau \textless 1/L\left(\nabla f \right)$
		\FOR{$ k=1,2,3,\ldots ,$MAXITER}
		\STATE $ \hat{\mathbf{x}}_{k}={\arg\mathop{\min} \limits_{\mathbf{y}}} \, \frac{1}{2\sigma^2}F_1\left( \mathbf{y} \right)-\nabla \frac{1}{2\sigma^2}F_2\left(\mathbf{y}_{k}\right)^T\left(\mathbf{y}-\mathbf{y}_{k}\right) $
		\STATE $ \mathbf{x}_{k} =\rm{ prox}_g\left(\hat{\mathbf{x}}_{k},\tau \right)
		= {\arg\mathop{\min} \limits_{\mathbf{x}}} \, \frac{\lambda \tau}{2} \mathbf{x}^T\left(\mathbf{x}-D\left(\mathbf{x} \right)\right)+\frac{1}{2}\left\| \mathbf{x}-\hat{\mathbf{x}}_{k}\right\|_2^2 $.  
		\STATE $ \alpha_{k+1} = \left( 1+\sqrt{1+4{\left(\alpha_{k}\right)}^2} \right)/2 $
		\STATE $ \mathbf{y}_{k+1} = \mathbf{x}_{k} + \frac{\alpha_{k}-1}{\alpha_{k+1}}\left(\mathbf{x}_{k} - \mathbf{x}_{k-1}\right) $
		\ENDFOR
		\STATE Output current point $\mathbf{x}_{MAXITER}$ as the approximate optimal solution.
	\end{algorithmic}
\end{algorithm}

\begin{algorithm}[H]
	\caption{prDeep-L2}
	\label{prDeep-L2}
	\begin{algorithmic}[1]
		\REQUIRE $\mathbf{y}_{1} = \mathbf{x}_{0} \in R^N, \alpha_{1}=1, \tau \textless 1/L\left(\nabla f \right)$
		\FOR{$ k=1,2,3,\ldots ,$MAXITER}
		\STATE $ \hat{\mathbf{x}}_{k}= \mathbf{y}_{k}-\frac{\tau}{2{\sigma}^2}\nabla F\left(\mathbf{y}_{k}\right) $
		\STATE $ 
		\mathbf{x}_{k} = { \rm{prox}_g}(\hat{\mathbf{x}}_{k},\tau) ={\arg\mathop{\min} \limits_{\mathbf{x}}} \, \frac{\lambda \tau}{2} \left\| \mathbf{x}-D\left(\mathbf{x}\right) \right\|_2^2+\frac{1}{2}\left\| \mathbf{x}-\hat{\mathbf{x}}_{k}\right\|_2^2 $. 
		\STATE $ \alpha_{k+1} = \left( 1+\sqrt{1+4{\left(\alpha_{k}\right)}^2} \right)/2 $
		\STATE $ \mathbf{y}_{k+1} = \mathbf{x}_{k} + \frac{\alpha_{k}-1}{\alpha_{k+1}}\left(\mathbf{x}_{k} - \mathbf{x}_{k-1}\right) $
		\ENDFOR
		\STATE Output current point $\mathbf{x}_{MAXITER}$ as the approximate optimal solution.
	\end{algorithmic}
\end{algorithm}

On the other hand, we choose the 2-norm regularization by denoising, which is $ \frac{\lambda}{2} \left\|\mathbf{x} -D\left(\mathbf{x} \right)\right\|_2^2$. Then the optimization problem can be written as
\begin{equation}
	{\arg\min_{\mathbf{x}}} \, \frac{1}{2{\sigma}^2} \left\| \mid \mathrm{\mathbf{A}} \mathbf{x} \mid - \mathbf{b}\right\|_2^2+\frac{\lambda}{2} \left\|\mathbf{x} -D\left(\mathbf{x} \right)\right\|_2^2. \label{the phase retrieval with and L_2}
\end{equation}

We use Algorithm \ref{prDeep-L2} to solve this optimization problem. The second step of Algorithm \ref{prDeep-L2} is
the same as the first step of Algorithm \ref{prDeep}. We use proximal operator for the third step of Algorithm \ref{prDeep-L2} such as $\mathbf{u}_1=\hat{\mathbf{x}}_{k}$, and then do the iterative process(26), where the last items in brackets are multiplication and division at corresponding positions respectively. The derivation process is as follows. The sub-optimization problem of the third step of Algorithm \ref{prDeep-L2} is that
\begin{equation}
	{\arg\mathop{\min} \limits_{\mathbf{x}}} \, \frac{\lambda \tau}{2} \left\| \mathbf{x}-D\left(\mathbf{x} \right) \right\|_2^2+\frac{1}{2}\left\| \mathbf{x}-\hat{\mathbf{x}}_{k}\right\|_2^2 .
\end{equation}
For convenience of description below, let 
\begin{equation}
	h\left(\mathbf{x} \right) = \frac{\lambda \tau}{2} \left\| \mathbf{x} -D\left(\mathbf{x} \right) \right\|_2^2+\frac{1}{2}\left\| \mathbf{x}-\hat{\mathbf{x}}_{k}\right\|_2^2 .
\end{equation}
Then the gradient of $h\left(\mathbf{x} \right)$ is obtained as following
\begin{equation}
	\begin{aligned}
		\nabla h\left(\mathbf{x} \right)=&\frac{\lambda \tau}{2}\left( 2\mathbf{x}-2D\left(\mathbf{x}\right) -2\mathbf{x} \cdot \mathbf{x} \cdot D\left(\mathbf{x} \right) +2\nabla D\left( \mathbf{x} \right) \cdot D\left( \mathbf{x} \right)  \right) \\
		&+ \mathbf{x}-\hat{\mathbf{x}}_{k} . \label{grad of h(x)} 
	\end{aligned}
\end{equation}

In\cite{romano2017little}, RED introduces the following property for denoiser D(x).
\begin{equation}
	\nabla D\left(\mathbf{x} \right) \cdot \mathbf{x} = D\left(\mathbf{x} \right) .
\end{equation}
Therefore, (\ref{grad of h(x)}) can be written as
\begin{equation}
	\nabla h\left( \mathbf{x} \right) = \left(\lambda \tau + 1\right) \mathbf{x} -\lambda \tau\left(2D\left( \mathbf{x} \right) - D\left( \mathbf{x} \right) \cdot D\left( \mathbf{x} \right) ./\mathbf{x} \right) - \hat{\mathbf{x}}_{k} .
\end{equation}
Let $\nabla h\left( \mathbf{x} \right) = 0$, we can get the following iterative steps.
\begin{equation}
	\mathbf{u}_j=\frac{1}{\lambda \tau + 1} \left( \hat{\mathbf{x}}_{k} + \lambda\tau \left(2D\left( \mathbf{u}_{j-1} \right) - D\left( \mathbf{u}_{j-1} \right) \cdot D\left( \mathbf{u}_{j-1} \right)./{\mathbf{u}_{j-1}}\right) \right) . \label{the inner loop}
\end{equation}

There is close relationship between the RED and 2-norm regularizations such that 
\begin{equation}
	\left\| \mathbf{x}-D\left( \mathbf{x} \right)\right\|^2 = \mathbf{x}^T\left(\mathbf{x} -D\left(\mathbf{x}\right)\right)+D\left(\mathbf{x} \right)^T\left(D\left(\mathbf{x} \right)-\mathbf{x}\right) .
\end{equation}

During the later numerical experiments on the 2-norm regularization by denoising, we firstly use once iteration with $\mathbf{x}^T\left(\mathbf{x}-D\left(\mathbf{x} \right)\right)$ as the predictive estimate of the 2-norm version. Secondly, we use prDeep-L2(\ref{prDeep-L2}) to
solve the sub-optimization problem with the 2-norm regularization by denoising. It is found that only one iteration can achieve excellent results.

\begin{table}[t]
	\centering
	\renewcommand{\tablename}{Table}
	\caption{The average PSNR/ average SSIM and average run time of the 128×128 test dataset UNT-6 under different noise levels.}
	\label{128_UNT-6_comparison}
	\begin{tabular}{ccccccc}
		\hline  
		{Method}  & \multicolumn{2}{c}{$\alpha=2$} & \multicolumn{2}{c}{$\alpha=3$} & \multicolumn{2}{c}{$\alpha=4$} \\ 
		{} & {PSNR/SSIM} & {TIME} & {PSNR/SSIM} & {TIME} & {PSNR/SSIM} & {TIME} \\ \hline
		{HIO} & 19.93/0.478 & 9.64 & 17.37/0.433 & 10.42 & 15.55/0.392 & 10.88 \\
		{prDeep} & 31.27/0.803 & 84.59 & 28.78/0.732 & 89.66 & 26.10/0.629 & 92.82 \\
		{prDeep-DC} & 32.21/0.829 & 204.76 & 29.14/0.736 & 205.54 & 27.44/0.676 & 220.11 \\	
		{prDeep-L2} & 32.26/0.827 & 177.50 & 29.25/0.749 & 159.67 & 27.43/0.698 & 183.07 \\ \hline
	\end{tabular}
\end{table}

\begin{table}[t]
	\centering
	\renewcommand{\tablename}{Table}
	\caption{The average PSNR/ average SSIM and average run time of the 128×128 test dataset NT-6 under different noise levels.}
	\label{128_NT-6_comparison}
	\begin{tabular}{ccccccc}
		\hline  
		{Method}  & \multicolumn{2}{c}{$\alpha=2$} & \multicolumn{2}{c}{$\alpha=3$} & \multicolumn{2}{c}{$\alpha=4$} \\ 
		{} & {PSNR/SSIM} & {TIME} & {PSNR/SSIM} & {TIME} & {PSNR/SSIM} & {TIME} \\ \hline
		{HIO} & 22.33/0.560 & 9.77 & 20.94/0.479 & 10.11 & 18.51/0.354 & 11.02 \\
		{prDeep} & 31.90/0.888 & 85.39 & 29.55/0.825 & 90.26 & 26.81/0.750 & 93.72 \\
		{prDeep-DC} & 32.54/0.900 & 197.70 & 30.20/0.846 & 211.17 & 28.00/0.775 & 212.35 \\ 
		{prDeep-L2} & 33.37/0.911 & 168.69 & 29.87/0.834 & 188.63 & 27.24/0.757 & 159.43 \\ \hline
	\end{tabular}
\end{table}

\begin{figure}[htbp]
	\centering
	\begin{subfigure}{0.19\linewidth}
		\caption*{Ground Truth}
		\centering
		\includegraphics[width=\linewidth]{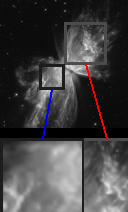}
		\caption*{PSNR/SSIM}
	\end{subfigure}
	\centering
	\begin{subfigure}{0.19\linewidth}
		\caption*{HIO}
		\centering
		\includegraphics[width=\linewidth]{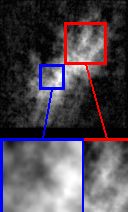}
		\caption*{21.17/0.388}
		
	\end{subfigure}
	\centering
	\begin{subfigure}{0.19\linewidth}
		\caption*{prDeep}
		\centering
		\includegraphics[width=\linewidth]{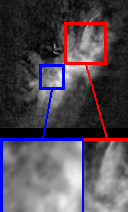}
		\caption*{21.18/0.408}
		
	\end{subfigure}
	\centering
	\begin{subfigure}{0.19\linewidth}
		\caption*{prDeep-DC}
		\centering
		\includegraphics[width=\linewidth]{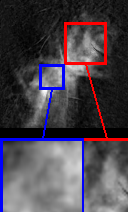}
		\caption*{24.75/0.560}
		
	\end{subfigure}
	\centering
	\begin{subfigure}{0.19\linewidth}
		\caption*{prDeep-L2}
		\centering
		\includegraphics[width=\linewidth]{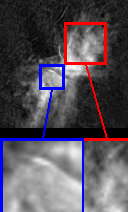}
		\caption*{24.21/0.543}
		
	\end{subfigure}
	
	\centering
	\begin{subfigure}{0.19\linewidth}
		\centering
		\includegraphics[width=\linewidth]{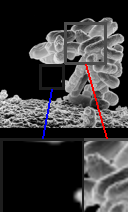}
		\caption*{PSNR/SSIM}
	\end{subfigure}
	\centering
	\begin{subfigure}{0.19\linewidth}
		\centering
		\includegraphics[width=\linewidth]{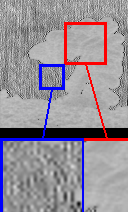}
		\caption*{5.91/0.340}
	\end{subfigure}
	\centering
	\begin{subfigure}{0.19\linewidth}
		\centering
		\includegraphics[width=\linewidth]{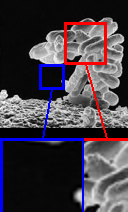}
		\caption*{37.50/0.980}
	\end{subfigure}
	\centering
	\begin{subfigure}{0.19\linewidth}
		\centering
		\includegraphics[width=\linewidth]{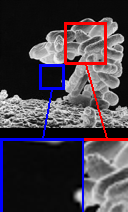}
		\caption*{37.63/0.982}
	\end{subfigure}
	\centering
	\begin{subfigure}{0.19\linewidth}
		\centering
		\includegraphics[width=\linewidth]{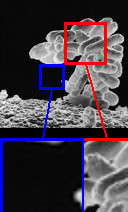}
		\caption*{37.64/0.969}
	\end{subfigure}
	
	\centering
	\begin{subfigure}{0.19\linewidth}
		\centering
		\includegraphics[width=\linewidth]{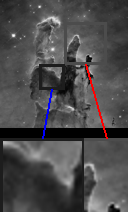}
		\caption*{PSNR/SSIM}
	\end{subfigure}
	\centering
	\begin{subfigure}{0.19\linewidth}
		\centering
		\includegraphics[width=\linewidth]{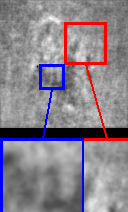}
		\caption*{23.14/0.440}
	\end{subfigure}
	\centering
	\begin{subfigure}{0.19\linewidth}
		\centering
		\includegraphics[width=\linewidth]{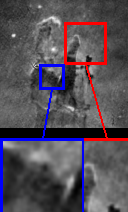}
		\caption*{27.72/0.729}
	\end{subfigure}
	\centering
	\begin{subfigure}{0.19\linewidth}
		\centering
		\includegraphics[width=\linewidth]{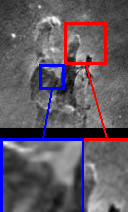}
		\caption*{28.14/0.721}
	\end{subfigure}
	\centering
	\begin{subfigure}{0.19\linewidth}
		\centering
		\includegraphics[width=\linewidth]{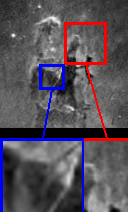}
		\caption*{28.25/0.719}
	\end{subfigure}
	
	\caption{Comparison of reconstruction of UNT-6 test images at $128 \times 128$ for Butterfly,Ecoli and Pillarsof at $\alpha=2$. }
	\label{128_alpha2_comparison Butterfly,Ecoli,Pillarsof}
\end{figure}

\begin{figure}[htbp]
	\centering
	\begin{subfigure}{0.19\linewidth}
		\caption*{Ground Truth}
		\centering
		\includegraphics[width=\linewidth]{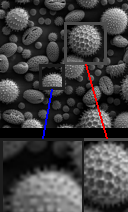}
		\caption*{PSNR/SSIM}
	\end{subfigure}
	\centering
	\begin{subfigure}{0.19\linewidth}
		\caption*{HIO}
		\centering
		\includegraphics[width=\linewidth]{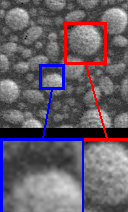}
		\caption*{23.80/0.718}
		
	\end{subfigure}
	\centering
	\begin{subfigure}{0.19\linewidth}
		\caption*{prDeep}
		\centering
		\includegraphics[width=\linewidth]{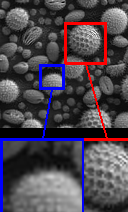}
		\caption*{29.12/0.912}
		
	\end{subfigure}
	\centering
	\begin{subfigure}{0.19\linewidth}
		\caption*{prDeep-DC}
		\centering
		\includegraphics[width=\linewidth]{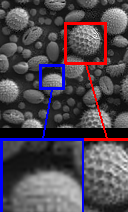}
		\caption*{29.64/0.920}
		
	\end{subfigure}
	\centering
	\begin{subfigure}{0.19\linewidth}
		\caption*{prDeep-L2}
		\centering
		\includegraphics[width=\linewidth]{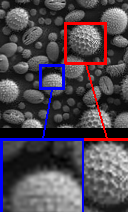}
		\caption*{30.66/0.938}
		
	\end{subfigure}
	
	\centering
	\begin{subfigure}{0.19\linewidth}
		\centering
		\includegraphics[width=\linewidth]{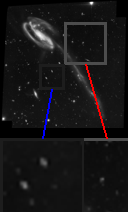}
		\caption*{PSNR/SSIM}
	\end{subfigure}
	\centering
	\begin{subfigure}{0.19\linewidth}
		\centering
		\includegraphics[width=\linewidth]{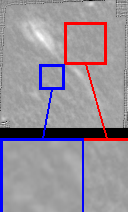}
		\caption*{22.07/0.464}
	\end{subfigure}
	\centering
	\begin{subfigure}{0.19\linewidth}
		\centering
		\includegraphics[width=\linewidth]{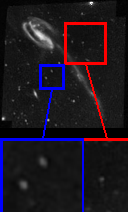}
		\caption*{35.45/0.906}
	\end{subfigure}
	\centering
	\begin{subfigure}{0.19\linewidth}
		\centering
		\includegraphics[width=\linewidth]{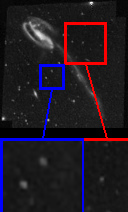}
		\caption*{36.26/0.907}
	\end{subfigure}
	\centering
	\begin{subfigure}{0.19\linewidth}
		\centering
		\includegraphics[width=\linewidth]{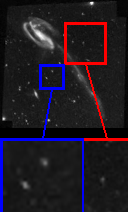}
		\caption*{36.35/0.918}
	\end{subfigure}
	
	\centering
	\begin{subfigure}{0.19\linewidth}
		\centering
		\includegraphics[width=\linewidth]{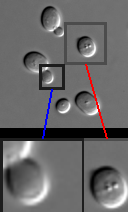}
		\caption*{PSNR/SSIM}
	\end{subfigure}
	\centering
	\begin{subfigure}{0.19\linewidth}
		\centering
		\includegraphics[width=\linewidth]{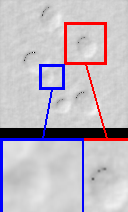}
		\caption*{23.49/0.519}
	\end{subfigure}
	\centering
	\begin{subfigure}{0.19\linewidth}
		\centering
		\includegraphics[width=\linewidth]{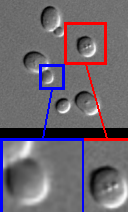}
		\caption*{36.63/0.884}
	\end{subfigure}
	\centering
	\begin{subfigure}{0.19\linewidth}
		\centering
		\includegraphics[width=\linewidth]{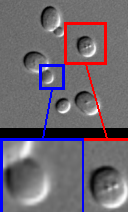}
		\caption*{36.84/0.885}
	\end{subfigure}
	\centering
	\begin{subfigure}{0.19\linewidth}
		\centering
		\includegraphics[width=\linewidth]{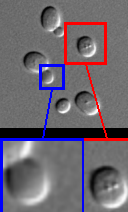}
		\caption*{36.48/0.877}
	\end{subfigure}
	
	\caption{Comparison of reconstruction of UNT-6 test images at $128 \times 128$ for Pollen,TadpoleGalaxy and Yeast at $\alpha=2$.}
	\label{128_alpha2_comparison Pollen,TadpoleGalaxy,Yeast}
\end{figure}

\begin{figure}[htbp]
	\centering
	\begin{subfigure}{0.19\linewidth}
		\caption*{Ground Truth}
		\centering
		\includegraphics[width=\linewidth]{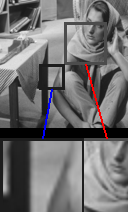}
		\caption*{PSNR/SSIM}

	\end{subfigure}
	\centering
	\begin{subfigure}{0.19\linewidth}
		\caption*{HIO}
		\centering
		\includegraphics[width=\linewidth]{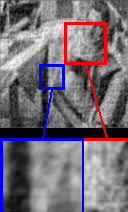}
		\caption*{22.99/0.665}
		
	\end{subfigure}
	\centering
	\begin{subfigure}{0.19\linewidth}
		\caption*{prDeep}
		\centering
		\includegraphics[width=\linewidth]{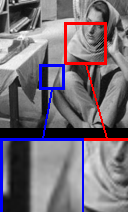}
		\caption*{35.59/0.974}
		
	\end{subfigure}
	\centering
	\begin{subfigure}{0.19\linewidth}
		\caption*{prDeep-DC}
		\centering
		\includegraphics[width=\linewidth]{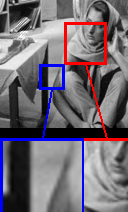}
		\caption*{35.66/0.976}
		
	\end{subfigure}
	\centering
	\begin{subfigure}{0.19\linewidth}
		\caption*{prDeep-L2}
		\centering
		\includegraphics[width=\linewidth]{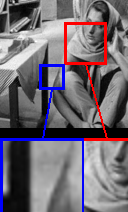}
		\caption*{37.89/0.980}
		
	\end{subfigure}
	
	\centering
	\begin{subfigure}{0.19\linewidth}
		\centering
		\includegraphics[width=\linewidth]{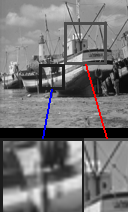}
		\caption*{PSNR/SSIM}
	\end{subfigure}
	\centering
	\begin{subfigure}{0.19\linewidth}
		\centering
		\includegraphics[width=\linewidth]{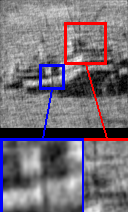}
		\caption*{22.39/0.508}
	\end{subfigure}
	\centering
	\begin{subfigure}{0.19\linewidth}
		\centering
		\includegraphics[width=\linewidth]{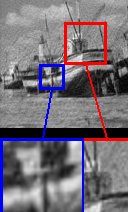}
		\caption*{27.87/0.803}
	\end{subfigure}
	\centering
	\begin{subfigure}{0.19\linewidth}
		\centering
		\includegraphics[width=\linewidth]{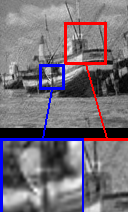}
		\caption*{28.04/0.812}
	\end{subfigure}
	\centering
	\begin{subfigure}{0.19\linewidth}
		\centering
		\includegraphics[width=\linewidth]{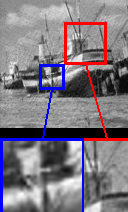}
		\caption*{29.51/0.851}
	\end{subfigure}
	
	\centering
	\begin{subfigure}{0.19\linewidth}
		\centering
		\includegraphics[width=\linewidth]{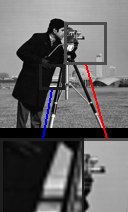}
		\caption*{PSNR/SSIM}
	\end{subfigure}
	\centering
	\begin{subfigure}{0.19\linewidth}
		\centering
		\includegraphics[width=\linewidth]{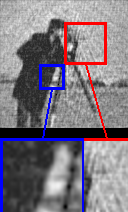}
		\caption*{21.87/0.459}
	\end{subfigure}
	\centering
	\begin{subfigure}{0.19\linewidth}
		\centering
		\includegraphics[width=\linewidth]{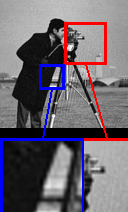}
		\caption*{34.73/0.924}
	\end{subfigure}
	\centering
	\begin{subfigure}{0.19\linewidth}
		\centering
		\includegraphics[width=\linewidth]{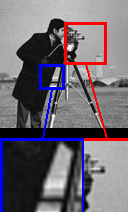}
		\caption*{37.60/0.949}
	\end{subfigure}
	\centering
	\begin{subfigure}{0.19\linewidth}
		\centering
		\includegraphics[width=\linewidth]{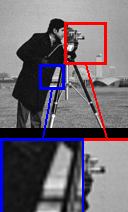}
		\caption*{37.81/0.954}
	\end{subfigure}
	
	\caption{Comparison of reconstruction of NT-6 test images at $128 \times 128$ for Barbara,Boat and Cameraman at $\alpha=2$.}
	\label{128_alpha2_comparison Barbara,Boat,Cameraman}
\end{figure}

\begin{figure}[htbp]
	\centering
	\begin{subfigure}{0.19\linewidth}
		\caption*{Ground Truth}
		\centering
		\includegraphics[width=\linewidth]{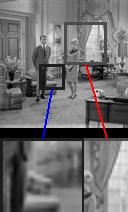}
		\caption*{PSNR/SSIM}

	\end{subfigure}
	\centering
	\begin{subfigure}{0.19\linewidth}
		\caption*{HIO}
		\centering
		\includegraphics[width=\linewidth]{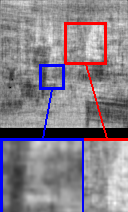}
		\caption*{20.95/0.478}
		
	\end{subfigure}
	\centering
	\begin{subfigure}{0.19\linewidth}
		\caption*{prDeep}
		\centering
		\includegraphics[width=\linewidth]{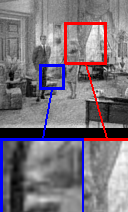}
		\caption*{29.46/0.863}
		
	\end{subfigure}
	\centering
	\begin{subfigure}{0.19\linewidth}
		\caption*{prDeep-DC}
		\centering
		\includegraphics[width=\linewidth]{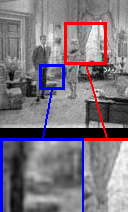}
		\caption*{29.52/0.867}
		
	\end{subfigure}
	\centering
	\begin{subfigure}{0.19\linewidth}
		\caption*{prDeep-L2}
		\centering
		\includegraphics[width=\linewidth]{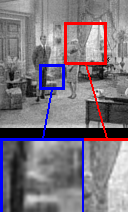}
		\caption*{30.46/0.886}
		
	\end{subfigure}
	
	\centering
	\begin{subfigure}{0.19\linewidth}
		\centering
		\includegraphics[width=\linewidth]{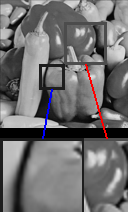}
		\caption*{PSNR/SSIM}
	\end{subfigure}
	\centering
	\begin{subfigure}{0.19\linewidth}
		\centering
		\includegraphics[width=\linewidth]{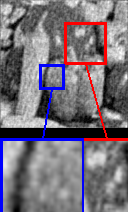}
		\caption*{23.22/0.634}
	\end{subfigure}
	\centering
	\begin{subfigure}{0.19\linewidth}
		\centering
		\includegraphics[width=\linewidth]{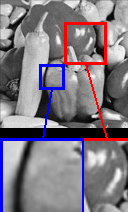}
		\caption*{37.69/0.970}
	\end{subfigure}
	\centering
	\begin{subfigure}{0.19\linewidth}
		\centering
		\includegraphics[width=\linewidth]{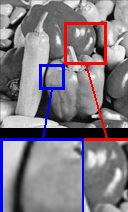}
		\caption*{37.70/0.969}
	\end{subfigure}
	\centering
	\begin{subfigure}{0.19\linewidth}
		\centering
		\includegraphics[width=\linewidth]{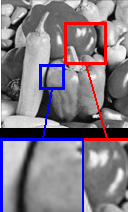}
		\caption*{37.48/0.967}
	\end{subfigure}
	
	\centering
	\begin{subfigure}{0.19\linewidth}
		\centering
		\includegraphics[width=\linewidth]{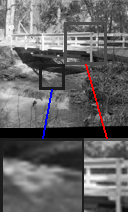}
		\caption*{PSNR/SSIM}
	\end{subfigure}
	\centering
	\begin{subfigure}{0.19\linewidth}
		\centering
		\includegraphics[width=\linewidth]{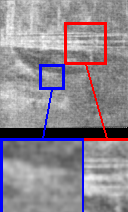}
		\caption*{22.57/0.618}
	\end{subfigure}
	\centering
	\begin{subfigure}{0.19\linewidth}
		\centering
		\includegraphics[width=\linewidth]{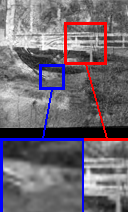}
		\caption*{26.06/0.798}
	\end{subfigure}
	\centering
	\begin{subfigure}{0.19\linewidth}
		\centering
		\includegraphics[width=\linewidth]{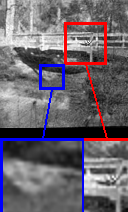}
		\caption*{26.72/0.824}
	\end{subfigure}
	\centering
	\begin{subfigure}{0.19\linewidth}
		\centering
		\includegraphics[width=\linewidth]{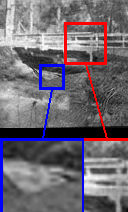}
		\caption*{27.06/0.826}
	\end{subfigure}
	
	\caption{Comparison of reconstruction of NT-6 test images at $128 \times 128$ for Couple,Peppers and Streamandbridge at $\alpha=2$.}
	\label{128_alpha2_comparison Couple,Peppers,Streamandbridge}
\end{figure}

\section{Numerical simulation}
In this section, we show the results of numerical experiments, introduce some experimental settings and show the restoration effects of different methods, including visual effects and quantitative values, to prove the effectiveness of our methods.

\subsection{Experimental setup}
\subsubsection{Measurement and Noise Model}
In this section, we test the algorithms with fourier measurement, which is widely used in many fields, including spectroscopy, optical microscope, remote sensing, chemical analysis and so on. At the same time, the main sources of niose are the poisson noise. We consider the following poisson noise model\cite{metzler2018prdeep} that
\begin{equation}
	\mathbf{y}^2 = \mid A \mathbf{x} \mid ^2+\mathbf{\omega} \; \rm{with} \; \mathbf{\omega} \sim N\left(0,\alpha ^2Diag\left(\mid A \mathbf{x} \mid ^2\right)\right) ,
\end{equation}
where ${\rm Diag}\left(\vert A\mathbf{x} \vert ^2\right)$ is the diagonal matrix of diagonal elements of $\vert A\mathbf{x} \vert ^2$, $	\mathbf{y}^2 / \alpha ^2 \sim {\rm Poisson}\left(\left(\vert A\mathbf{x}\vert / \alpha\right)^2\right)$, $\alpha$ is the parameter of Poisson noise, and $	\mathbf{y}^2$ is the intensity information with Poisson noise.

\subsubsection{Parameter tuning and initialization}
We use Fourier measurement with the oversampling rate of 4 and the test images are $128 \times 128$ grayscale images NT-6 and UNT-6. For the regularization term parameter, we choose $\sigma_\omega$, where $\sigma_\omega^2$ represents the variance of the sample noise. Because the fourier measurement method is sensitive to the selection of initial values, the HIO algorithm  is used to select an initial value. Our algorithms iterate 200 times, the sub-optimization problem using DC iterates 10 times, and the proximal gradient formula of g(x) iterates once. We refer to prDeep\cite{metzler2018prdeep} and  6 natural images and 6 unnatural images are used as test datasets, which are denoted by NT-6(Figure.\ref{NT-6}) and UNT-6(Figure.\ref{UNT-6}) respectively. They obey different distributions. We consider three different  possion noise levels with $ \alpha \in \left[ 2,3,4 \right] $. 

\subsection{Experimental analysis of phase retrieval based on prDeep-DC}
Table.\ref{128_UNT-6_comparison} and Table.\ref{128_NT-6_comparison} 
present the average PSNR, SSIM, and running time for $128 \times 128$ 
test datasets UNT-6 and NT-6, respectively, under different noise levels. HIO performs worst across all noise levels, likely due to its inability to solve noisy phase retrieval. prDeep-DC consistently outperforms prDeep in all aspects. For UNT-6, prDeep-DC yields PSNR improvements of 0.94dB, 0.36dB, and 1.34dB over prDeep. For NT-6, prDeep-DC shows PSNR gains of 0.64dB, 0.65dB, and 1.19dB. This demonstrates prDeep-DC’s more prominent reconstruction effect, particularly at higher noise levels.

Figures \ref{128_alpha2_comparison Butterfly,Ecoli,Pillarsof}, \ref{128_alpha2_comparison Pollen,TadpoleGalaxy,Yeast}, \ref{128_alpha2_comparison Barbara,Boat,Cameraman}, and \ref{128_alpha2_comparison Couple,Peppers,Streamandbridge} show the reconstruction comparisons for test images UNT-6 and NT-6, respectively, when $\alpha=2$. Visually, HIO provides the worst reconstruction quality. In contrast, prDeep, prDeep-DC, and prDeep-L2 all produce reconstructions very close to the ground truth. Importantly, prDeep-DC’s reconstruction maps retain finer details and are more faithful to the ground truth than those of prDeep.

\subsection{Experimental analysis of phase retrieval based on prDeep-L2}
An examination of the numerical results presented in Tables \ref{128_UNT-6_comparison} and \ref{128_NT-6_comparison}
reveals that prDeep-L2 consistently outperforms prDeep across all evaluated aspects. Specifically, for the UNT-6 test images, prDeep-L2 achieves superior PSNR values by 0.99dB, 0.47dB, and 1.33dB compared to prDeep under the three distinct noise levels, respectively. Similarly, for the NT-6 test images, prDeep-L2 demonstrates PSNR gains of 1.47dB, 0.32dB, and 0.43dB over prDeep across the same noise conditions. These quantitative findings indicate that prDeep-L2 yields a more pronounced reconstruction effect for UNT-6 test images, particularly at higher noise levels. Conversely, for NT-6 test images, prDeep-L2 exhibits its most prominent advantage at lower noise levels. The final column in Figures 
\ref{128_alpha2_comparison Butterfly,Ecoli,Pillarsof}, \ref{128_alpha2_comparison Pollen,TadpoleGalaxy,Yeast}, \ref{128_alpha2_comparison Barbara,Boat,Cameraman} and \ref{128_alpha2_comparison Couple,Peppers,Streamandbridge}
displays the restoration map generated by prDeep-L2. A visual inspection confirms that prDeep-L2’s reconstructed maps also preserve finer details and exhibit greater fidelity to the ground truth compared to those obtained by prDeep.

\section{Conclusion}
This paper addresses noisy phase retrieval by proposing a DC programming and DnCNN-based denoising regularization framework. First, the objective function is decomposed into a difference of two convex functions, leading to the novel prDeep-DC algorithm, which combines DCA and prDeep. Second, we introduce prDeep-L2 by applying a 2-norm denoising regularization term to phase retrieval. Numerical experiments confirm the excellent performance of our methods, both quantitatively and visually.

\section*{Acknowledgments}
This research work was supported in part by the National Natural Science Foundation of China (grants 12261037).
\bibliographystyle{unsrt}
\bibliography{Phase_Retrieval_Based_on_DC_and_DnCNN}




\end{document}